\documentclass[12pt]{article}
\usepackage{amsfonts}
\usepackage{psfrag}
\newtheorem{lemma}{\sc Lemma}[section]

\newtheorem{definition}[lemma]{\sc Definition}
\newtheorem{proposition}[lemma]{\sc Proposition}
\newtheorem{theorem}[lemma]{\sc Theorem}
\newtheorem{remark}[lemma]{\sc Remark}

\newtheorem{corollary}[lemma]{\sc Corollary}

\def\scaledpicture #1 by #2 (#3 scaled #4){{
  \dimen0=#1 \dimen1=#2
  \divide\dimen0 by 1000 \multiply\dimen0 by #4
  \divide\dimen1 by 1000 \multiply\dimen1 by #4
  \picture \dimen0 by \dimen1 (#3 scaled #4)}
  }

\def\CC{\mathbb{C}}

\def\SS{\mathbb{S}}
\def\RR{\mathbb{R}}

\def\NN{\mathbb{N}}

\def\bb{\mathtt{b}}

\def\dd{\mathtt{d}}
\def\ii{\mathtt{i}}

\def\pp{\mathtt{p}}

\def\ss{\mathtt{s}}

\def\rr{\mathtt{r}}

\def\qed{\vrule height5pt depth0pt width5pt} 
\def\[{[\![}
\def\]{]\!]}

\def\truc#1{\mathop{\triangleright}\limits_{\raise5pt\hbox{$_{#1}$}}}

\def\0{{\bf 0}}
\def\1{{\bf 1}}
\def\2{{\bf 2}}
\def\3{{\bf 3}}

\def\ddd{\dd}\def\rrr{\rr}

\begin{document}

\baselineskip=22pt

\centerline{\large\bf Infinite products of $2\times2$ matrices }

\centerline{\large\bf and the Gibbs properties of Bernoulli convolutions}
\bigskip

\bigskip

\centerline{by {\sc Eric Olivier \&
Alain Thomas} }
\bigskip

{\leftskip=2cm
\rightskip2cm

{\small
{\bf Abstract.--} {We consider the infinite sequences 
$(A_n)_{n\in\NN}$ of $2\times2$ matrices with nonnegative entries, where the $A_n$ are taken in a finite set of matrices. Given a vector
$V=\pmatrix{v_1\cr v_2}$ with $v_1,v_2>0$, we give a necessary and
sufficient condition for $\displaystyle{A_1\dots A_nV\over\vert\vert A_1\dots A_nV\vert\vert}$ to converge  uniformly. In
application we prove that the Bernoulli convolutions related to the numeration in Pisot quadratic bases are weak Gibbs.}\par
\par}\par
}

\bigskip


\centerline{\large \bf Introduction}
\medskip
Let ${\cal M}=\{M_0,\dots,M_{\ss-1}\}$ be a finite subset of the set -- stable by matrix multiplication -- of nonnegative
and column-allowable $d\times d$ matrices (i.e., the matrices with nonnegative entries and without null column). We associate to any sequence
$\left(\omega_n\right)_{n\in\NN}$ with terms in
${\cal S}:=\{0,1,\dots,\ss-1\}$, the sequence of product matrices
$$
P_n(\omega)=M_{\omega_1}M_{\omega_2}\dots M_{\omega_n}.
$$
Experimentally, in most cases each normalized column of $P_n(\omega)$ converges when $n\to~\infty$ to a limit-vector, which depends on $\omega\in{\cal
S}^{\NN}$ and may depend on the index of the column.

Nevertheless the normalized 
rows of
$P_n(\omega)$ in general do not converge: suppose for instance that all the matrices in $\cal M$ are positive but do not have the same positive normalized
left-eigenvector, let $L_k$ such that $L_kM_k=\rho_kL_k$. For any positive matrix $M$, the normalized rows of $M{M_0}^n$
converge to $L_0$ and the ones of
$M{M_1}^n$ to $L_1$. Consequently we can choose the sequence $(n_k)_{k\in\NN}$ sufficiently increasing such that the normalized rows of
$M_0^{n_1}M_1^{n_2}\dots M_0^{n_{2k-1}}$ converge to
$L_0$ while the ones of
$M_0^{n_1}M_1^{n_2}\dots M_0^{n_{2k-1}}M_1^{n_{2k}}$ converge to~$L_1$. This proves -- if $L_0\ne L_1$ -- that the normalized rows in
$P_n(\omega)$ do not converge when $\omega=0^{n_1}1^{n_2}0^{n_3}1^{n_4}\dots $. 

Now in case $\cal M$ is a set of positive matrices it is clear that, if both normalized columns and normalized rows in
$P_n(\omega)$ converge then -- after
replacing each matrix
$M_k$ by $\displaystyle{1\over\rho_k}M_k$~-- the matrix $P_n(\omega)$ itself converges: the previous counterexample proves that
the matrices $P_n(\omega)$ have a common left-eigenvector for any $n$, and a straightforward computation (using 
the limits of the normalized columns in $P_n(\omega)$) proves the existence of $\displaystyle\lim_{n\to\infty}P_n(\omega)$.

The existence of a common left-eigenvector is settled in a more general context by L.~Elsner and S. Friedland 
(\cite[Theorem 1]{EF}), in case $\cal M$ is a finite set of matrices with entries in~$\CC$. This theorem means (after transposition of the matrices) that if 
$P_n(\omega)$ converges to a non-null limit, then there exists $N\in\NN$ such that the matrices $M_{\omega_n}$ for $n\ge N$ have a common 
left-eigenvector for the eigenvalue $1$. Now, L. Elsner \& S. Friedland (in \cite[Main Theorem]{EF}) and I. Daubechies \& J. C. Lagarias (in \cite[Theorem 5.1]{DL'}
\hbox{(resp. \cite[Theorem~4.2]{DL})}) give necessary and sufficient conditions for $P_n(\omega)$ to converge for any
$\omega\in{\cal S}^\NN$ (resp., to converge to a continuous map).

By these theorems we see that the problem of the convergence of the normalized columns in $P_n(\omega)$ is very different from the problem of the convergence of
$P_n(\omega)$ itself. Let for instance
$M_0=\pmatrix{1/2&1/2\cr1/3&2/3}$ and
$M_1=\pmatrix{1/2&1/2\cr1/2&1/2}$; then the normalized columns in
$$
P_n(\omega)=\left\{\begin{array}{ll}\displaystyle{1\over10}\cdot\pmatrix{4+6\cdot6^{-n}
&6-6\cdot6^{-n}\cr4-4\cdot6^{-n}
&6+4\cdot6^{-n}}&\hbox{if }\omega_1\dots\omega_n=0\dots0\\
\displaystyle{1\over10}\cdot\pmatrix{4+6^{-h}
&6-6^{-h}\cr4+6^{-h}
&6-6^{-h}}&\hbox{if }\omega_1\dots\omega_n=\omega_1\dots\omega_{n-h-1}10\dots0\end{array}\right.
$$
converge to $\pmatrix{1/2\cr1/2}$ for any $\omega\in\{0,1\}^{\NN}$, but $P_n(\omega)$ diverges (altough it is bounded) if $\omega$~is not eventually constant. 

In Section 1 we study the uniform convergence -- in direction -- of $P_n(\omega)V$ in case the $M_k$ are $2\times 2$ nonnegative
column-allowable  matrices and $V=\pmatrix{v_1\cr v_2}$ a positive vector (Theorem \ref{convergence}). Notice that the convergence in direction of the columns of $P_n(\omega)$, to a same vector, implies the ones of $P_n(\omega)V$, but the converse is not true: see for instance the case
${\cal M}=\left\{\pmatrix{2&0\cr1&1}\right\}$.

The second section is devoted to the {\it Bernoulli convolutions} \cite{E}, which have been studied
since the early 1930's (see \cite{PSS} for
the other references). We give a matricial relation for such measures.

In the third section we apply more precisely Theorem \ref{convergence} to prove that certain Bernoulli convolutions are weak Gibbs in the following sense
(see \cite{Y}): given a system off affine contractions $\SS_\varepsilon:\RR\to\RR$ such that the intervals
$\SS_\varepsilon([0,1])$ make a partition of $[0,1]$ for~$\varepsilon\in{\cal S}=\{0,1,\dots,\ss-1\}$, a measure $\eta$ supported by $[0,1]$ is weak Gibbs
w.r.t. $\{\SS_\varepsilon\}_{\varepsilon=0}^{\ss-1}$ if there exists a
map $\Phi:{\cal S}^\NN\to\RR$, continuous for the product topology, such that
\begin{equation}\label{A2}
\lim_{n\to\infty}\left({\eta\[\xi_1\dots\xi_n\]\over
\exp\Big(\sum_{k=0}^{n-1} \Phi(\sigma^k\xi)\Big)}\right)^{1/n}=1\quad\hbox{uniformly on }\xi\in{\cal S}^\NN,
\end{equation}
where $\[\xi_1\dots\xi_n\]:=\SS_{\xi_1}\circ\dots\circ\SS_{\xi_n}([0,1])$ and $\sigma$ is the shift on ${\cal S}^\NN$.
Let us give a sufficient condition for $\eta$ to be weak Gibbs.
For each $\xi\in{\cal S}^\NN$ we put
$\phi_1(\xi)=\log\eta\[\xi_1\]$ and for~$n\ge2$,
\begin{equation}\label{nstep}
\phi_n(\xi)=
\log\left({\eta\[\xi_1\cdots\xi_n\]\over\eta\[\xi_2\cdots\xi_n\]}\right).
\end{equation}
The continuous map $\phi_n:{\cal S}^\NN\to\RR$ ($n\ge1$) is 
the {\em $n$-step potential} of $\eta$. Assume the existence of the uniform limit $\Phi=\lim_{n\to\infty}\phi_n$; it is then straightforward that for~$n\ge1$,
\begin{equation}\label{A*}
{1\over K_n}\leq{\eta\[\xi_1\cdots\xi_n\]\over
\exp\Big(\sum_{k=0}^{n-1} \Phi(\sigma^k\xi)\Big)}\leq K_n 
\quad\hbox{with}\quad K_n=\exp\left(\sum_{k=1}^n\Vert\Phi-\phi_n\Vert_\infty\right).
\end{equation}
By a well known lemma on the Ces\`aro sums, $K_1,K_2,\dots$ form a subexponential
sequence of positive real numbers, that is
$\lim_{n\to\infty}\left(K_n\right)^{1/n}=1$ and thus, (\ref{A*}) means 
$\eta$ is
weak Gibbs  w.r.t. $\{\SS_\varepsilon\}_{\varepsilon=0}^{\ss-1}$.

Now the weak Gibbs property can be proved for certain Bernoulli
convolutions by computing the $n$-step potential by means of products of matrices (see \cite{FO} for the Bernoulli convolution associated with the golden
ratio $\displaystyle\beta={1+\sqrt5\over2}$ -- called the {\em Erd\"os measure} -- and the application to the multifractal analysis). In Theorem \ref{quadratic} we
generalize this result in case $\beta>1$ is a quadratic number with conjugate $\beta'\in]-1,0[$.

\begin{section}{Infinite product of $2\times2$ matrices}

From now the vectors $X=\pmatrix{x_1\cr x_2}$ and the matrices $A=\pmatrix{a&b\cr c&d}$ we consider are supposed to be nonnegative and column-allowable that is,
$x_1,x_2,a,b,c,d$ are nonnegative and $x_1+x_2,a+c,b+d$ are positive. In particular we suppose that the matrices in ${\cal M}=\{M_0,\dots,M_{\ss-1}\}$ satisfy these
conditions. We associate to $X$ the normalized vector:
$$
\mathtt{N}(X):=\pmatrix{{x_1\over x_1+x_2}\cr{x_2\over x_1+x_2}}=\pmatrix{\mathtt{n}(X)\cr1-\mathtt{n}(X)}\quad\hbox{where}\quad\mathtt{n}(X):={x_1\over
x_1+x_2}$$
and define the distance between the column of $A$ (or the rows of ${^tA}$):
$$d_{\rm{columns}}(A):=\left\vert\mathtt{n}\left(\pmatrix{a\cr c}\right)-\mathtt{n}\left(\pmatrix{b\cr
d}\right)\right\vert ={\vert\det A\vert\over(a+c)(b+d)}=:d_{\rm{rows}}({^tA}).
$$

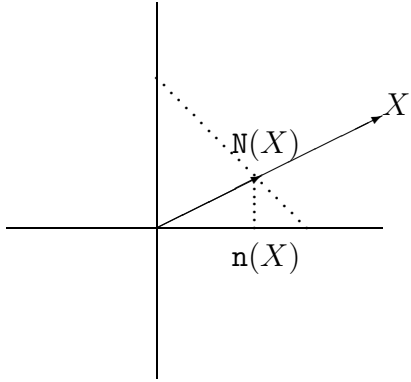
\setlength{\unitlength}{1cm}
\begin{picture}(6,6.5)(5,4)
\put(5,7){\line(1,0){5}}
\put(7,5){\line(0,1){5}}
\put(10,8.5){$X$}
\put(8,8){$\mathtt{N}(X)$}
\put(8,6.5){$\mathtt{n}(X)$}
\put(7,7){\vector(2,1){3}}
\put(7,7){\vector(2,1){1.4}}
\put(7,9){\circle*{0}}
\put(7.1,8.9){\circle*{0}}
\put(7.2,8.8){\circle*{0}}
\put(7.3,8.7){\circle*{0}}
\put(7.4,8.6){\circle*{0}}
\put(7.5,8.5){\circle*{0}}
\put(7.6,8.4){\circle*{0}}
\put(7.7,8.3){\circle*{0}}
\put(7.8,8.2){\circle*{0}}
\put(7.9,8.1){\circle*{0}}
\put(8,8){\circle*{0}}
\put(8.1,7.9){\circle*{0}}
\put(8.2,7.8){\circle*{0}}
\put(8.3,7.7){\circle*{0}}
\put(8.4,7.6){\circle*{0}}
\put(8.5,7.5){\circle*{0}}
\put(8.6,7.4){\circle*{0}}
\put(8.7,7.3){\circle*{0}}
\put(8.8,7.2){\circle*{0}}
\put(8.9,7.1){\circle*{0}}
\put(9,7){\circle*{0}}
\put(8.3,7.525){\circle*{0}}
\put(8.3,7.42){\circle*{0}}
\put(8.3,7.315){\circle*{0}}
\put(8.3,7.21){\circle*{0}}
\put(8.3,7.105){\circle*{0}}
\put(8.3,7){\circle*{0}}
\end{picture}

\begin{theorem}\label{convergence}Given $V=\pmatrix{v_1\cr v_2}$ with $v_1,v_2>0$, the sequence of vectors $\mathtt{N}(P_n(\omega)V)$ converges uniformly for
$\omega\in{\cal S}^{\NN}$ only in the five following cases:

\underline{Case 1}: $\pmatrix{a&b\cr0&d}\in{\cal M}\Rightarrow
a\ge d;\enskip\pmatrix{a&0\cr c&d}\in{\cal M}\Rightarrow a\le d;$

$\pmatrix{a&b\cr c&0},\pmatrix{0&b'\cr c'&d'}\in{\cal
M}\Rightarrow bc'\ge b'c$;\enskip
no matrix in ${\cal
M}$ has the form
$\pmatrix{a&0\cr 0&d}$~or~$\pmatrix{0&b\cr c&0}$.

\underline{Case 2}: $\pmatrix{a&b\cr0&d}\in{\cal M}\Rightarrow
a<d;\enskip\pmatrix{a&0\cr c&d}\in{\cal M}\Rightarrow a>d;$

$\pmatrix{a&b\cr c&0},\pmatrix{0&b'\cr c'&d'}\in{\cal
M}\Rightarrow bc'<b'c$.

\underline{Case 3}: $\pmatrix{a&b\cr0&d}\in{\cal M}\Rightarrow
a\ge d;\enskip\pmatrix{a&0\cr c&d}\in{\cal M}\Rightarrow a>d$;

no matrix in ${\cal
M}$ has the form~$\pmatrix{0&b\cr c&d}$.

\underline{Case 4}: $\pmatrix{a&b\cr0&d}\in{\cal M}\Rightarrow
a<d;\enskip\pmatrix{a&0\cr c&d}\in{\cal M}\Rightarrow a\le d$;

no matrix in ${\cal
M}$ has the form~$\pmatrix{a&b\cr c&0}$;

\underline{Case 5}: $V$ is an eigenvector of all the matrices in $\cal M$.

\end{theorem}

\begin{corollary}\label{corol0}
If $\mathtt{N}(P_n(\cdot)V)$ converges uniformly on ${\cal
S}^{\NN}$, the limit
do not depend on the positive vector $V$, except in the fifth case of Theorem \ref{convergence}.
\end{corollary}
{\it Proof.} Suppose that ${\cal M}$ satisfies the conditions of the case 1,2,3 or 4 in Theorem \ref{convergence} and let $V,W$ be two positive vectors. Then the following set ${\cal M}'$ also do:

${\cal M}':={\cal M}\cup\{M_s\}$, where $M_s$ is the matrix
whose both columns are $W$.

Denoting by $\omega'=\omega_1\dots\omega_n\overline s$ the sequence defined by\quad
$\omega'_i=\cases{\omega_i&if $i\le n$\cr s&if $i> n$\cr}$
one has for~any~$\omega\in{\cal S}^{\NN}$
$$\mathtt{N}(P_n(\omega)V)-\mathtt{N}(P_n(\omega)W)=\mathtt{N}(P_n(\omega')V)-\mathtt{N}(P_{n+1}(\omega')V)
$$
and this tends to $0$, according to the uniform Cauchy property of the sequence
$\mathtt{N}(P_n(\cdot)V)$.~\qed

Nevertheless, this limit may depend of $V$ if one assume only that $V$ is nonnegative. For instance, if ${\cal
M}=\left\{\pmatrix{1&1\cr0&2},\pmatrix{1&1\cr1&1}\right\}$ then
$\displaystyle\lim_{n\to\infty}\mathtt{N}\left(P_n(\omega)\pmatrix{1\cr1}\right)=\pmatrix{1/2\cr1/2}$ differs from
$\lim_{n\to\infty}\displaystyle\mathtt{N}\left(P_n(\omega)\pmatrix{1\cr0}\right)$ iff $\omega=\overline 0$ (implying the second limit is not uniform on ${\cal
S}^{\NN}$).

\begin{subsection}{Geometric considerations}
We follow the ideas of E. Seneta about products of nonnegative matrices in Section 3 of~\cite{S}, or stochastic matrices in Section 4. 
In what follows we denote the matrices by $A=~\pmatrix{a&b\cr c&d}$, $A'=~\pmatrix{a'&b'\cr c'&d'}$ or
$A_n=~\pmatrix{a_n&b_n\cr c_n&d_n}$ for
$n\in\NN$, and we suppose they are nonnegative and column-allowable. We define the coefficient
$$
\tau(A):=\sup_{d_{\rm{columns}}(A')\ne0}{d_{\rm{columns}}(A'A)\over d_{\rm{columns}}(A')}.
$$
The straightforward formula
\begin{equation}\label{columndist}
d_{\rm{columns}}\left(\prod_{k=1}^n(A_k)\right)\le d_{\rm{columns}}(A_1)\prod_{k=2}^n\tau(A_k)
\end{equation}
is of use to prove Theorem \ref{convergence} because, according to the following proposition one has $\tau(A)<1$ if $A$ is positive.
\begin{proposition}\label{computau}
$$
\tau(A)=\left\{\begin{array}{ll}\displaystyle{\left\vert\sqrt{ad}-\sqrt{bc}\right\vert\over\sqrt{ad}+\sqrt{bc}}&\hbox{if $A$ do not have any null
row}\\&\\0&\hbox{otherwise.}\end{array}\right.
$$
\end{proposition}
{\it Proof.} $\displaystyle{d_{\rm{columns}}(A'A)\over d_{\rm{columns}}(A')}={\vert\det
A\vert\over(a+c/x)(bx+d)}$$\left(\hbox{where }\displaystyle x={a'+c'\over b'+d'}\right)$ is maximal for
$\displaystyle x=\sqrt{cd\over ab}$.~\qed
\begin{remark}One can consider -- instead of $d_{\rm{columns}}$ -- the angle between the columns of~$A$:
$$
\displaystyle\alpha(A):=
\left\vert\arctan{a\over c}-\arctan{b\over d}\right\vert,
$$
or the Hilbert distance between the columns of a positive matrice $A$:
$$
d_{\rm{Hilbert}}(A):=
\left\vert\log{a\over c}-\log{b\over d}\right\vert.
$$
This last can be interpreted  either as
the distance between the columns or the rows of $A$, because $d_{\rm{Hilbert}}(A)=d_{\rm{Hilbert}}({^tA})$. The Birkhoff coefficient
$\displaystyle\tau_{\rm{Birkhoff}}(A):=\sup_{d_{\rm{Hilbert}}(A')\ne0}{d_{\rm{Hilbert}}(A'A)\over d_{\rm{Hilbert}}(A')}
$ has -- from \cite[Theorem~(3.12)]{S} -- the same value as $\tau(A)$ in Proposition \ref{computau}, and probably as a large class of coefficients defined in this
way.
\end{remark}

In the following proposition we list the properties of $d_{\rm{columns}}$ that are required for proving Theorem \ref{convergence}.
\begin{proposition}\label{list}(i) $\displaystyle\sup_{d_{\rm{columns}}(A')\ne0}{d_{\rm{columns}}(AA')\over
d_{\rm{columns}}(A')}={\vert\det A\vert\over\min\left((a+c)^2,(b+d)^2\right)}=:\tau_1(A)$.

(ii) If $A$ is positive then $\displaystyle\sup_{d_{\rm{columns}}(A')\ne0}{d_{\rm{rows}}(AA')\over
d_{\rm{columns}}(A')}\le{\vert\det A\vert\over\min(a,b)\cdot\min(c,d)}=:\tau_2(A)$.

(iii) If $\displaystyle\lim_{n\to\infty}d_{\rm{columns}}(A_n)=0$ then
$\displaystyle\lim_{n\to\infty}d_{\rm{columns}}(AA_nA')=0$ and, assuming that $A$ is positive, $\displaystyle\lim_{n\to\infty}d_{\rm{rows}}(AA_nA')=0$.

(iv) Suppose the matrices $A_n$ are upper-triangular.
If $\displaystyle\inf_{k\in{\NN}}{a_k\over d_k}\ge~1$ and
$\displaystyle\sum_{k\in{\NN}}{b_k\over d_k}=\infty$ then
$$\lim_{n\to\infty}d_{\rm{columns}}(A_1\dots A_n)=
\lim_{n\to\infty}d_{\rm{columns}}(A_n\dots A_1)=
0.$$

(v) Suppose the matrices $A_n$ are lower-triangular. If $\displaystyle\inf_{k\in{\NN}}{d_k\over a_k}\ge~1$ and
$\displaystyle\sum_{k\in{\NN}}{c_k\over a_k}=\infty$ then
$$\lim_{n\to\infty}d_{\rm{columns}}(A_1\dots A_n)=
\lim_{n\to\infty}d_{\rm{columns}}(A_n\dots A_1)=
0.$$

\end{proposition}
{\it Proof.} (i) and (ii) are obtained from the formula
$$
d_{\rm{columns}}(AA')={\det A\cdot\det A'\over((a+c)a'+(b+d)c')\cdot((a+c)b'+(b+d)d')},
$$
and the relation $d_{\rm{rows}}(AA')=d_{\rm{columns}}({^tA'}\ {^tA})$.

(iii) is due to the fact that the inequalities of items (i), (ii) and (\ref{columndist}) imply $d_{\rm{columns}}(AA_nA')\le\tau_1(A)d_{\rm{columns}}(A_n)\tau(A')$ and -- if $A$ is positive -- $d_{\rm{rows}}(AA_nA')\le\tau_2(A)d_{\rm{columns}}(A_n)\tau(A')$.

(iv) follows from the formula
$$A_1\dots A_n=\pmatrix{a_1\dots a_n&s_n\cr0&d_1\dots d_n},$$where $\displaystyle
s_n=\sum_{k=1}^na_1\dots a_{k-1}b_kd_{k+1}\dots d_n\ge d_1\dots d_n\sum_{k=1}^n{b_k\over d_k}$.

(v) can be deduced from (iv) by using the relation $\pmatrix{0&1\cr1&0}\pmatrix{a&0\cr c&d}\pmatrix{0&1\cr1&0}=\pmatrix{d&c\cr0&a}$.~\qed

We need also the following:
\begin{proposition}\label{eigen}Let $V_A$ be a nonnegative eigenvector associated to the maximal eigenvalue of $A$, and $C$ a cone of nonnegative
vectors containing $V_A$. If $\det A\ge0$ then $C$ is stable by left-multiplication by $A$.
\end{proposition}
{\it Proof.} The discriminant of the characteristic polynomial of $A$ is $(a-d)^2+4bc$. In case this discriminant is null the proof is obtained by direct
computation, because $A=\pmatrix{a&b\cr 0&a}$ or $\pmatrix{a&0\cr c&a}$. Otherwise $A$ has two eigenvalues $\lambda>\lambda'$ and, given a nonnegative
vector $X$, there exists a real $\alpha$
and an eigenvector
$W_A$ (associated to~$\lambda'$) such that
$$
X=\alpha V_A+W_A\quad\hbox{and}\quad AX=\lambda\alpha V_A+\lambda' W_A=\lambda'X+(\lambda-\lambda')\alpha V_A.
$$
Notice that $\alpha\ge0$ (because the nonnegative vector $A^nX=\lambda^n\alpha V_A+{\lambda' }^nW_A$ converges in direction to $\alpha V_A$) 
and $\lambda'\ge0$ (from the hypothesis $\det A\ge0$). Hence $AX$ is a nonnegative linear combination of $X$ and $V_A$; if $X$ belongs to $C$ then $AX$ also
do.~\qed

\end{subsection}

\begin{subsection}{How pointwise convergence implies uniform convergence}
Let $m$ and $M$ be the bounds of $\mathtt{n}(P_n(\omega)V)$ for $n\in\NN$ and $\omega\in{\cal S}^{\NN}$, and let $M_V:=\pmatrix{m&M\cr1-m&1-M}$.
Each real $x\in[m,M]$ can be written $x=mx_1+Mx_2$ with $x_1,x_2\ge0$ and $x_1+x_2=1$; in particular the real $x=\mathtt{n}(P_n(\omega)V)$ can be written in this
form, hence
\begin{equation}\label{c™ne}
\forall\omega\in{\cal S}^{\NN},\ \exists t_1,t_2\ge0,\ P_n(\omega)V=M_V\pmatrix{t_1\cr t_2}.
\end{equation}

\begin{proposition}\label{loc}
If $d_{\rm{columns}}(P_n(\cdot)M_V)$ converges pointwise to $0$ when $n\to\infty$, then $\mathtt{N}(P_n(\cdot)V)$
converges uniformly on ${\cal S}^{\NN}$.
\end{proposition}
{\it Proof.} Suppose the pointwise convergence holds. Given $\omega\in{\cal S}^{\NN}$ and $\varepsilon>0$, there exists the integer
$n=n(\omega,\varepsilon)$ such that
$d_{\rm{columns}}(P_n(\omega)M_V)\le\varepsilon$.
The family of cylinders
$C(\omega,\varepsilon):=\[\omega_1\dots\omega_{n(\omega,\varepsilon)}\]$, for
$\omega$ running over ${\cal S}^{\NN}$, is a covering of the compact~${\cal S}^{\NN}$; hence there exists a finite subset $X\subset{\cal S}^{\NN}$ such that 
${\cal S}^{\NN}=\bigcup_{\omega\in X}C(\omega,\varepsilon)$. Let $p>q\ge n_\varepsilon:=\max\{n(\omega,\varepsilon)\;;\;\omega\in~X\}$. 
For any $\xi\in{\cal S}^{\NN}$,
there exists $\zeta\in X$ such that $\xi\in C(\zeta,\varepsilon)$ that is, $\xi_k=\zeta_k$ for any $k\le n=n(\zeta,\varepsilon)$. 
From (\ref{c™ne}) there exists two nonnegative vectors $V_p$ and $V_q$ such that $P_p(\xi)V=P_n(\zeta)M_VV_p$ and $P_q(\xi)V=P_n(\zeta)M_VV_q$. Denoting by $M(p,q)$
the column-allowable matrix whose columns are $V_p$ and $V_q$ we have -- in view of (\ref{columndist})
$$
\begin{array}{rcl}\left\vert\mathtt{n}\left(P_p(\xi)V\right)-\mathtt{n}\left(P_q(\xi)V\right)\right\vert
&=&d_{\rm{columns}}\left(P_n(\zeta)M_VM(p,q)\right)\\&\le&d_{\rm{columns}}\left(P_n(\zeta)M_V\right)\\&\le&\varepsilon,\end{array}
$$
implying the uniform Cauchy property for $\mathtt{N}(P_n(\cdot)V)$.~\qed

\end{subsection}

\begin{subsection}{Proof of the uniform convergence of $\mathtt{N}(P_n(\cdot)V)$}
According to Proposition~\ref{loc} it is sufficient to prove that
$\displaystyle\lim_{n\to\infty}d_{\rm{columns}}(P_n(\omega)M_V)=0$ for each
$\omega\in{\cal S}^{\NN}$. This convergence is obvious in the following cases:

$\bullet$ If there exists $N$ such that $M_{\omega_N}$ has rank $1$, then $P_n(\omega)M_V$ has rank $1$ for $n\ge N$ and 
$$
\forall n\ge N,\ d_{\rm{columns}}(P_n(\omega)M_V)=0.
$$

$\bullet$ If there exists infinitely many integers $n$ such that $M_{\omega_n}$ is a positive matrix, one has
$\displaystyle\tau(M_{\omega_n})\le\rho:=\max_{M\in{\cal M},\ M>0}\tau(M)<1$, and the formula (\ref{columndist})
implies
$$
\lim_{n\to\infty}d_{\rm{columns}}(P_n(\omega)M_V)=0.
$$

$\bullet$ Similarly, this limit is null also in case there exists infinitely many integers $n$ such that $M_{\omega_n}M_{\omega_{n+1}}$ is a positive matrix.

So we can make from now the following hypotheses on the sequence $\omega$ under consideration:

(H): $\det M_{\omega_n}\ne0$ for any $n\in\NN$, and there exists $N$ such that the matrix $M_{\omega_n}M_{\omega_{n+1}}$ 
has at least one null entrie for any $n>N$.

Proof in the case 1: Since the couples of matrices $\pmatrix{a&b\cr c&0},\pmatrix{0&b'\cr c'&d'}\in{\cal
M}$ satisfy $\displaystyle{b\over c}\ge{b'\over c'}$, there exists a real $\alpha$ such that 
$$
\forall\pmatrix{a&b\cr c&0},\pmatrix{0&b'\cr c'&d'}\in{\cal
M},\quad{b\over c}\ge\alpha\ge{b'\over c'}.
$$
Let $\Delta=\pmatrix{0&\alpha\cr1&0}$. We denote by $\cal P$ the set of $2\times2$ matrices with
nonnegative determinant and by $\tilde{\cal M}$ the subset of $\cal P$ defined as follows:
$$
\tilde{\cal M}:=\{\Delta^{-1}M,\ M\Delta\ ;\ M\in{\cal M}\setminus{\cal P}\}\cup\{M,\ \Delta^{-1}M\Delta\ ;\ M\in{\cal M}\cap{\cal P}\}.
$$
This set of matrices also satisfies the
conditions mentionned in the case 1: for instance if $\pmatrix{a&b\cr c&0}\in{\cal M}$, the matrix $\Delta^{-1}\pmatrix{a&b\cr c&0}=\pmatrix{c&0\cr
a/\alpha&b/\alpha}$ satisfies $c\le b/\alpha$, and so one. For any sequence $\varepsilon_0,\varepsilon_1,\varepsilon_2,\dots$ of elements of $\{0,1\}$ such that
$\varepsilon_0=0$ we can write

\begin{equation}\label{Delta}\begin{array}{rcl}
P_n(\omega)&=&M_{\omega_1}M_{\omega_2}\dots
M_{\omega_n}\\&=&\left(\Delta^{-\varepsilon_0}M_{\omega_1}\Delta^{\varepsilon_1}\right)\cdot\left(\Delta^{-\varepsilon_1}M_{\omega_2}\Delta^{\varepsilon_2}\right)
\cdot\dots
\cdot\left(\Delta^{-\varepsilon_{n-1}}M_{\omega_n}\Delta^{\varepsilon_n}\right)\cdot\Delta^{-\varepsilon_n}
\\&=&A_1A_2\dots A_n\Delta^{-\varepsilon_n}\end{array}
\end{equation}

where $A_n:=\Delta^{-\varepsilon_{n-1}}M_{\omega_n}\Delta^{\varepsilon_n}$ for any $n\in\NN$.
By the following choice of the sequence $(\varepsilon_n)_{n\in\NN}$, the matrices
$A_n$ belong to $\tilde{\cal M}$:
$$
\varepsilon_n=\left\{\begin{array}{ll}\varepsilon_{n-1}&\hbox{if }\det M_{\omega_n}>0\\1-\varepsilon_{n-1}&\hbox{otherwise.}\end{array}\right.
$$

The hypotheses (H) imply that either all the matrices $A_n$ for $n>N$ are upper-triangular, or all of them are lower-triangular (otherwise 
$M_{\omega_n}M_{\omega_{n+1}}=\Delta^{\varepsilon_{n-1}}A_nA_{n+1}\Delta^{-\varepsilon_{n+1}}$ is positive for some $n>N$). By Proposition~\ref{list}~(iv)~and~(v),
$$\lim_{n\to\infty}d_{\rm{columns}}(A_{N+1}\dots A_{n})=0.$$
From (\ref{Delta}) and Proposition~\ref{list} (iii),
$\displaystyle\lim_{n\to\infty}d_{\rm{columns}}(P_n(\omega)M_V)=~0$.

Proof in the case 2: We use the matrix $\Delta$ and the set of matrices $\tilde{\cal M}$ defined in the previous case; here the real $\alpha$ 
is supposed such that
$\displaystyle{b\over c}\le\alpha\le{b'\over c'}$ for any $\pmatrix{a&b\cr c&0},\pmatrix{0&b'\cr c'&d'}\in{\cal
M}$, and consequently $\tilde{\cal M}$ satisfies the hypotheses of the case 2.
This imply that each matrix in $\tilde{\cal M}$ has a positive eigenvector. Let $C$ be the (minimal) cone containing $V$, $\Delta^{-1}V$ and the
positive eigenvectors of the matrices in $\tilde{\cal M}$. From (\ref{Delta}) and Proposition \ref{eigen}, $P_n(\omega)V$ belongs to this cone for any
$\omega\in{\cal S}^{\NN}$ hence $M_V$ is positive.

Using again the relation (\ref{Delta}) we have
\begin{equation}\label{rowedP_n}
d_{\rm{columns}}\left(P_n(\omega)M_V\right)=d_{\rm{rows}}\left({^tM_V}\ {^t\Delta^{-\varepsilon_n}}\ {^tA_n}\dots{^tA_1}\right).
\end{equation}
Each matrix ${^tA_n}$ for $n>N$ satisfy $a>d$ 
if ${^tA_n}$ is upper-triangular, and
$a<d$ if it is lower-triangular. By Proposition~\ref{list}~(iv)~and~(v),
$\lim_{n\to\infty}d_{\rm{columns}}({^tA_n}\dots{^tA_{N+1}})=0$. This implies $\displaystyle\lim_{n\to\infty}d_{\rm{columns}}(P_n(\omega)M_V)=~0$
by applying Proposition \ref{list} (iii) to the r.h.s. in~(\ref{rowedP_n}).

Proof in the case 3: Let $C'$ be the (minimal) cone containing $V$, the nonnegative eigenvectors (associated to the maximal eigenvalues) of the
matrices in ${\cal M}\cap{\cal P}$, and the column-vectors of the matrices in ${\cal M}\setminus{\cal P}$. All the vectors delimiting $C'$ are distinct from
$\pmatrix{0\cr1}$, and Proposition
\ref{eigen} implies that
$P_n(\omega)V\in C'$ for any
$\omega\in{\cal S}^{\NN}$. Hence $m$ and $M$ that is, the bounds of $\mathtt{n}(P_n(\omega)V)$), are positive.

Suppose first that $M_{\omega_n}$ is lower-triangular for any $n\in\NN$ and let $\pmatrix{\alpha_n&0\cr\gamma_n&\delta_n}=P_n(\omega)$. The hypotheses
of the case 3 imply $\displaystyle\lim_{n\to\infty}{\delta_n\over\alpha_n}=0$. A simple computation gives
$\displaystyle d_{\rm{columns}}\left(P_n(\omega)M_V\right)\le{\delta_n\over\alpha_n}\cdot{M-m\over Mm}$ hence
$\displaystyle\lim_{n\to\infty}d_{\rm{columns}}\left(P_n(\omega)M_V\right)=0$. This conclusion remains valid if $M_{\omega_n}$ is eventually lower-triangular. 

Suppose now $M_{\omega_n}$ is not lower-triangular for infinitely many $n$. The hypotheses mentionned in the case 3 and (H) imply that $M_{\omega_n}$ is
upper-triangular for any $n>N$ (because for each $n$ such that $M_{\omega_n}$ is lower-triangular or has the form
$\pmatrix{a&b\cr c&0}$, (H) implies that $M_{\omega_{n+1}}$ is lower-triangular). Proposition~\ref{list}~(iii)~and(iv) implies that
$\displaystyle\lim_{n\to\infty}d_{\rm{columns}}(P_n(\omega)M_V)=~0$.

Proof in the case 4: Let ${\cal M}'$ be the set of matrices $M'_k=\Delta^{-1}M_k\Delta$ for $k=0,\dots,\ss-1$, and
let $V'=\Delta^{-1}V$ (here we can choose $\Delta=\pmatrix{0&1\cr1&0}$). The set ${\cal M}'$ satisfies the hypotheses of the case 3 hence
$\displaystyle\lim_{n\to\infty}d_{\rm{columns}}(P_n(\omega)M_V)=\lim_{n\to\infty}d_{\rm{columns}}(\Delta M'_{\omega_1}\dots M'_{\omega_n}V')=~0$.

\end{subsection}

\begin{subsection}{Proof of the converse assertion in Theorem \ref{convergence}}
Now we suppose the existence of the uniform limit $\displaystyle V(\cdot):=\lim_{n\to\infty}\mathtt{N}(P_n(\cdot)V)$ and we want to check the conditions contained
in one of the five cases involved in Theorem \ref{convergence}.  Let ${\cal M}^2$ be the set of matrices $MM'$ for $M,M'\in{\cal
M}$, and let
$\cal U$ (resp.
$\cal L$) be the set of upper-triangular (resp. lower-triangular) matrices $M\in{\cal M}\cup{\cal M}^2$.

We first prove that $\cal U$ cannot contain a couple of matrices $A=\pmatrix{a&b\cr 0&d}$ and $A'=~\pmatrix{a'&b'\cr 0&d'}$ such that $a\ge d$ and $a'<d'$: 
suppose that $\cal U$ contain such matrices let, for simplicity, $M_0=A$ and $M_1=A'$.
One has $V(\overline0)=\displaystyle\lim_{n\to\infty}\mathtt{N}(A^nV)$, and this limit is also the normalized nonnegative right-eigenvector of $A$ associated to its
maximal eigenvalue, hence
$V(\overline0)=\pmatrix{1\cr0}$. Similarly, $V(\overline1)$ is colinear to $\pmatrix{b'\cr d'-a'}$ (eigenvector of $A'$) hence distinct from $V(\overline0)$.
Moreover, for fixed $N\in\NN$
\begin{equation}\label{discontinuity}\begin{array}{rcl}
\displaystyle V(1^N\overline0)=\lim_{n\to\infty}\mathtt{N}(A'^NA^nV)&=&\displaystyle\lim_{n\to\infty}\mathtt{N}(A'^N\mathtt{N}(A^nV))\\
&=&\mathtt{N}\left(A'^NV(\overline0)\right)\\&=&V(\overline0).
\end{array}\end{equation}
Since $1^N\overline0$ tends to $\overline1$ when $N\to\infty$, the inequality $V(\overline0)\ne V(\overline1)$ contradicts the continuity of the map $V$.
This proves that the couple of matrices
$A,A'\in{\cal U}$ such that
$a\ge d$ and
$a'<d'$ do not exist. Similarly, the couple of matrices $A,A'\in\cal L$ such that $a\le d$ and $a'>d'$ do not exist.

$\bullet$ Suppose that all the matrices in $\cal U$ satisfy $a\ge d$ and all the ones in $\cal L$ satisfy $a\le~d$. If
$\cal M$ contains a matrix of the form $\pmatrix{a&0\cr0&d}$, it is necessarily an homothetic matrix. If it contains a matrix of the form $\pmatrix{0&b\cr c&0}$,
the square of this matrix is homothetic. So in both cases ${\cal M}\cup{\cal M}^2$ contains an homothetic matrix, let $H$. We use the same method as above: since the map $V$ is continuous, the vector
$\displaystyle\lim_{n\to\infty}\mathtt{N}(H^nV)$ must be equal~to $\displaystyle\lim_{N\to\infty}\left(\lim_{n\to\infty}\mathtt{N}(H^NMH^nV)\right)$ for any $M\in{\cal M}$. But the first is $\mathtt{N}(V)$ and the second $\mathtt{N}(MV)$, hence $V$ is an eigenvector of all the matrices in~$\cal M$. Suppose now that $\cal M$ do not contain matrices of the form $\pmatrix{a&0\cr0&d}$ nor $\pmatrix{0&b\cr c&0}$: all the conditions of the case 1
are satisfied.

$\bullet$ Suppose that all the matrices in $\cal U$ satisfy $a<d$ and all the ones in $\cal L$ satisfy $a>d$; then the conditions of the case 2
are satisfied.

$\bullet$ Suppose that all the matrices $A=\pmatrix{a&b\cr0&d}\in{\cal U}$ satisfy $a\ge d$ and all the matrices $A'=\pmatrix{a'&0\cr c'&d'}\in{\cal L}$ satisfy
$a'>d'$. If there exists $A\in{\cal U}$, $A'\in{\cal L}$ and $M=\pmatrix{0&\beta\cr\gamma&\delta}\in{\cal M}$, the map
$V$ is discontinuous because
$$
\lim_{n\to\infty}\mathtt{N}\left(A'^NMA^n\pmatrix{v_1\cr v_2}\right)=\lim_{n\to\infty}
\mathtt{N}\left(\pmatrix{a'&0\cr c'&d'}^N\pmatrix{\beta&0\cr\delta&\gamma}\pmatrix{d&0\cr b&a}^n\pmatrix{v_2\cr v_1}\right)=\pmatrix{0\cr1}
$$
differs from $\displaystyle\lim_{n\to\infty}\mathtt{N}\left(A'^n\pmatrix{v_1\cr v_2}\right)$ which is colinear to $\pmatrix{a'-d'\cr c'}$. Hence, either $\cal M$ do
not contain a matrix of the form
$\pmatrix{0&\beta\cr\gamma&\delta}$ and we are in the case 3, or ${\cal U}=\emptyset$ and we are in the case 2, or ${\cal L}=\emptyset$ and we are in the case 1.

$\bullet$ The case when all the matrices $\pmatrix{a&b\cr0&d}\in\cal U$ satisfy $a<d$ and all the matrices $\pmatrix{a'&0\cr c'&d'}\in\cal L$ satisfy $a'\le d'$ 
is symmetrical to the previous, by using the set of matrices ${\cal M}':=\{\Delta^{-1}M\Delta\ ;\ M\in{\cal M}\}$ and
the vector $V'=\Delta^{-1}V$, where $\Delta=\pmatrix{0&1\cr1&0}$.

\end{subsection}

\end{section}

\begin{section}{Some properties of the Bernoulli convolutions in base $\beta>1$}
Given a real $\beta>1$, an integer $\dd>~\beta$ and a $\dd$-dimensional
probability vector $\pp:=(\pp_i)_{i=0}^{\ddd-1}$, the
{\em $\pp$-distributed} \hbox{\em $(\beta,\dd)$-Bernoulli convolution} is
by definition the probability distribution $\mu_\pp$
of the random variable $X$ defined by
$$\forall\omega\in{\cal D}^{\NN}:=~\{\0,\dots,\dd-1\}^{\NN},\ X(\omega)=\sum\nolimits_{k=1}^\infty{\omega_k\over\beta^k},$$
where $\omega\mapsto\omega_k$ ($k=1,2,\cdots$) is a sequence of i.i.d. random
variables assuming the values
$i=\0,\1,\dots,\dd-1$ with probability $\pp_i$.

Denoting by
$\overline\omega$ the sequence such that
$\overline\omega_k=\dd-1-\omega_k$ for any $k$, one has the relation $\displaystyle X(\omega)+X(\overline\omega)=\alpha:={\dd-1\over\beta-1}$. Hence, setting 
$\overline\pp_i=\pp_{\ddd-1-i}$ for any $i=\0,\1,\dots,\dd-1$, the following
symmetry relation holds for any Borel set
$B\subset\RR$:
\begin{equation}\label{sym}
\mu_\pp(B)=\mu_{\overline\pp}(\alpha-B)
\end{equation}
(notice that the support of
$\mu_\pp$ is a subset of $[0,\alpha]$).

The measure $\mu_\pp$ also satisfy the following selfsimilarity relation: denoting by $\sigma$ the~shift~on~${\cal D}^\NN$ one has -- for any Borel set
$B\subset\RR$
$$
X(\omega)\in B\Leftrightarrow X(\sigma\omega)\in(\beta B-\omega_0)
$$
hence, using the independance of the random variables $\omega\mapsto\omega_k$,
\begin{equation}\label{A0}
\mu_\pp(B)=\sum_{k=0}^{\ddd-1}\pp_k\cdot
\mu_\pp(\beta B-k)\quad\hbox{for any Borel set $B\subset\RR$}
\end{equation}
and in particular
\begin{equation}\label{A1}
\mu_\pp(B)=\pp_0\cdot\mu_\pp(\beta B)\quad\hbox{if }\beta B\subset[0,1].
\end{equation}
The following proposition is
proved in \cite[Theorem~2.1 and Proposition~5.4]{DST} in case the probability vector $\pp$ is uniform:
\begin{proposition}
The $1$-periodic map $H:]-\infty,0]\to\RR$ defined by
$$
H(x)=(\pp_0)^x\cdot\mu_\pp\left(\left[0,\beta^x\right[\right)
$$
is continuous and a.e. differentiable. Moreover $H$ is not differentiable on a certain continuum of points if $\beta$ is an irrational Pisot number or an
integer and -- in this latter case~-- if $\beta$ do not divide
$\dd$.
\end{proposition}

Let us give also the matricial form of the relation (\ref{A0}) (from \cite[\S 2.1]{OST}).
We define the (finite or countable) set ${\cal I}_{(\beta,\ddd)}=\{0=\ii_0,\ii_1,\cdots\}$ as follows (where $\cal
B$ is the alphabet $\{0,1,\dots,\bb-1\}$ such that $\bb-1<\beta\le\bb$):
%

\begin{definition}\label{DefI}
{\sl
${\cal I}_{(\beta,\ddd)}$ is the set of
$i\in]-1,\alpha[$ for which there exists
$-1<i_1,\cdots,i_n<\alpha$ with
$0\triangleright i_1\triangleright\cdots
\triangleright i_n\triangleright i$,
where $x\triangleright y$ means
that exists $(\varepsilon,k)\in{\cal
B}\times{\cal D}$ such that $y=\beta x +(\varepsilon-k)$.}
\end{definition}

%
Let $\varepsilon\in{\cal B}$; the entries of the matrix $M_\varepsilon$
are -- for the row index $i$ and the column index $j$, with $\ii_i,\ii_j\in{\cal I}_{(\beta,\ddd)}$,
$$
M_\varepsilon(i,j)=
\cases{
\pp_k & if $k=\varepsilon+\beta \ii_i-\ii_j\in{\cal D}$\cr
0&otherwise.\cr}
$$
%

Setting $\displaystyle\RR_\varepsilon(x)={x+\varepsilon\over\beta}$
for any $\varepsilon\in{\cal B}$ and $x\in\RR$, we have the following

\begin{proposition}(\cite[Lemma 2.2]{OST})
\label{Ath4}
{\sl If ${\cal I}_{(\beta,\ddd)}=
\big\{\ii_0,\cdots,\ii_{\rrr-1}\big\}$
then, for any Borel set $B\subset[0,1]$ and any $\varepsilon\in{\cal B}$ such that ${\RR_\varepsilon}^{-1}(B)\subset[0,1]$,
$$
\pmatrix{
\mu_\pp(B+\ii_0)\cr
\vdots\cr
\mu_\pp(B+\ii_{\rrr-1})}
=M_\varepsilon
\pmatrix{
\mu_\pp\Big({\RR_\varepsilon}^{-1}(B)+\ii_0\Big)\cr
\vdots\cr
\mu_\pp\Big({\RR_\varepsilon}^{-1}(B)+\ii_{\rrr-1}\Big)\cr}.
$$}
\end{proposition}

\begin{remark}The finiteness of ${\cal I}_{(\beta,\ddd)}$ is assured, according to \cite[\S 2.2]{OST}, if $\beta$ is an irrational Pisot number or an
integer.\end{remark}
We shall use also the probability distribution
of the fractionnal part of the random variable $X$, that we denote by $\mu_\pp^*$. Suppose that $\alpha$ belongs to $]1,2[$, or equivalently
that $\beta<\dd<2\beta-1$. Then $\mu_\pp^*$ -- which has support
$[0,1]$ -- satisfy the following relation for any Borel set
$B\subset[0,1]$:
$$
\mu_\pp^*(B)=\mu_\pp(B)+\mu_\pp(B+1)
$$
and, if $B\subset[\alpha-1,1]$,
\begin{equation}\label{muŽtoile}
\mu_\pp^*(B)=\mu_\pp(B).
\end{equation}

The following proposition points out that in certain cases, the restriction of
$\mu_\pp$ (or $\mu_\pp^*$) to the interval $[\alpha-1,1]$ is ''representative'' of
$\mu_\pp$ itself.
\begin{proposition}\label{reduce}
Suppose $\displaystyle\beta<\dd\le\beta+1-{1\over\beta}$.

(i) The interval $]0,\alpha[$ is the reunion of $I_k:=\displaystyle\left[{1\over\beta^{k+1}}\ ,\ {1\over\beta^k}\right]$
and
$I'_k:=\displaystyle\left[\alpha-{1\over\beta^k}\ ,\ \alpha-{1\over\beta^{k+1}}\right]$ for $k\in\NN\displaystyle\cup\{0\}$

(ii) Let $B\subset\RR$ be a Borel
set . If $B\subset I_k$ (or equivalently if $\alpha-B\subset I'_k$), then $\beta^kB$ and $\alpha-\beta^kB$ are two subsets of $[\alpha-1,1]$ such that
$$
\begin{array}{rcl}
\mu_\pp(B)&=&\pp_0^k\cdot\mu_\pp^*(\beta^kB)\\
\mu_\pp(\alpha-B)&=&\pp_{\dd-1}^k\cdot\mu_\pp^*(\alpha-\beta^kB).
\end{array}
$$
\end{proposition}
{\it Proof.} (i) The hypothesis on $\dd$ implies $\alpha<2$ hence $]0,\alpha[$ is the reunion of $]0,1]$ and $[\alpha-1,\alpha[$.

(ii) $B\subset I_k\Rightarrow\beta^kB\subset\displaystyle\left[{1\over\beta}\ ,\ 1\right]\subset[\alpha-1,1]$. The equality
$\mu_\pp(B)=\pp_0^k\cdot\mu_\pp^*(\beta^kB)$ results from (\ref{A1}) and (\ref{muŽtoile}).

Since $\beta^kB\subset[\alpha-1,1]$ one has $\alpha-\beta^kB\subset[\alpha-1,1]$. The equality $\mu_\pp(\alpha-B)=\pp_{\dd-1}^k\cdot~\mu_\pp(\alpha-~\beta^kB)$
follows from (\ref{sym}), (\ref{A1}) and (\ref{muŽtoile}).~\qed

\end{section}

\begin{section}{Bernoulli convolution in Pisot quadratic bases}

In this section $\beta>1$ is solution of the equation
$x^2=ax+b$ (with integral $a$ and $b$), and we suppose that the
other solution belongs to ]-1,0[.
This implies $\displaystyle1\le b\le a\le\beta-{1\over\beta}<\beta<a+1$.
Let $\pp=(\pp_0,\dots,\pp_a)$ be a positive probability vector; the Bernoulli convolution $\mu_\pp$ has support $[0,\alpha]$, where
$\displaystyle\alpha={a\over\beta-1}$ belongs to $]1,2[$. The condition in Proposition \ref{reduce} is satisfied hence
it is sufficient to study the Gibbs properties of $\mu_\pp^*$ on its support $[0,1]$, to get the local properties of $\mu_\pp$ on
$[0,\alpha]$ (see \cite{FO} for the multifractal analysis of the weak Gibbs measures).

With the notations of the previous subsection one has ${\cal B}={\cal D}=\{0,\dots,a\}$,
\hfill\break${\cal I}_{(\beta,a+1)}=\{0,1,\beta-a\}$ and -- for any $\varepsilon\in{\cal B}$
$$
M_\varepsilon=
\pmatrix{\pp_\varepsilon&\pp_{\varepsilon-1}&0\cr0&0&\pp_{a+\varepsilon}\cr \pp_{b+\varepsilon}&\pp_{b+\varepsilon-1}&0},
$$
where, by convention, $\pp_i=0$ if $i\not\in{\cal D}$.

Notice that the intervals $\RR_\varepsilon([0,1])$ do not make a partition of $[0,1]$ for $\varepsilon\in{\cal B}$ but, setting
$$
\SS_\varepsilon:=\cases{\RR_\varepsilon&for $0\le\varepsilon\le a-1$\cr\RR_a\circ\RR_{\varepsilon-a}&for $a\le\varepsilon\le
a+b-1$}
$$
the intervals $\SS_\varepsilon([0,1])$ make such a partition for $\varepsilon\in{\cal A}:=\big\{0,\cdots,a+b-1\big\}$.

\begin{theorem}\label{quadratic}{\sl The measure $\mu_\pp^*$ is weak Gibbs w.r.t.
$\{\SS_\varepsilon\}_{\varepsilon=0}^{a+b-1}$ if~and~only~if
$\pp_0^2\ge~{\pp_a}\ \pp_{b-1}$ and $\pp_0\ \pp_{a-b+1}\ge~{\pp_a}^2$.}
\end{theorem}
{\it Proof.} 
The $n$-step potential of $\mu_\pp^*$ can be computed by means of the matrices 
$$
M^*_\varepsilon:=\cases{\displaystyle M_\varepsilon&for $0\le\varepsilon\le a-1$\cr\displaystyle M_aM_{\varepsilon-a}&for $a\le\varepsilon\le
a+b-1$.}
$$
Indeed by applying Proposition \ref{Ath4} to the sets $B=\[\xi_1\dots\xi_n\]$ and $B'=\[\xi_2\dots\xi_n\]$, one~has
\begin{equation}\label{label}
\exp(\phi_n(\xi))=\log\left({\pmatrix{1&1&0}M^*_{\xi_1}\dots M^*_{\xi_n}V\over\pmatrix{1&1&0}M^*_{\xi_2}\dots M^*_{\xi_n}V}\right),
\quad\hbox{where}\quad V:=
\pmatrix{
\mu_\pp([0,1])\cr
\mu_\pp([0,1]+1)\cr
\mu_\pp([0,1]+\beta-a)\cr}.
\end{equation}

Now the matrices $M^*_\varepsilon$ are $3\times3$ and we shall use $2\times2$ ones. 
The matrices defined -- for any $\varepsilon\in{\cal A}':=\{0,\dots,2a\}$ -- by
$$
M^\star_\varepsilon:=\left\{\begin{array}{ll}M_0M_\varepsilon&\hbox{if }\varepsilon\le a\\M_{\varepsilon-a}&\hbox{if }\varepsilon>a
\end{array}\right.
$$
satisfy the commutation relation $YM^\star_\varepsilon=P_\varepsilon Y$, where
$$
Y:=\pmatrix{1&0&0\cr0&1&0}\quad\hbox{and}\quad
P_\varepsilon:=\left\{\begin{array}{ll}\pmatrix{\pp_0\pp_\varepsilon&\pp_0\pp_{\varepsilon-1}\cr\pp_a\pp_{b+\varepsilon}&\pp_a\pp_{b+\varepsilon-1}\cr}&\hbox{if }\varepsilon\le a\\
\pmatrix{\pp_{\varepsilon-a}&\pp_{\varepsilon-a-1}\cr0&0\cr}&\hbox{if }\varepsilon>a.
\end{array}\right.
$$
Let $\xi\in{\cal A}^\NN$ such that $\sigma\xi\ne\overline0$. There exists an integer $k\ge0$ and 
$\varepsilon\in{\cal A}\setminus\{0\}$ such that
$$
M_{\xi_2}\dots M_{\xi_{k+2}}={M_0}^kM_\varepsilon.
$$
One can associate to the sequence $\xi$, the sequence $\zeta\in{{\cal A}'}^\NN$ such that

$$
\forall n\ge k+4,\ \exists k(n)\in\NN,\quad M^*_{\xi_{k+3}}\dots M^*_{\xi_n}=M^\star_{\zeta_1}\dots
M^\star_{\zeta_{k(n)}}\quad\hbox{or}\quad M^\star_{\zeta_1}\dots
M^\star_{\zeta_{k(n)}}M_0.
$$
Now -- according to (\ref{label}) and the commutation relation
\begin{equation}\label{over}
n\ge k+4\Rightarrow\exp(\phi_n(\xi))={\pmatrix{1&1&0}M^*_{\xi_1}M_0^kQ_\varepsilon\mathtt{N}(P_{\zeta_1}\dots
P_{\zeta_{k(n)}}W)\over\pmatrix{1&1&0}M_0^kQ_\varepsilon\mathtt{N}(P_{\zeta_1}\dots
P_{\zeta_{k(n)}}W)}
\end{equation}
where $Q_\varepsilon:=\pmatrix{\pp_\varepsilon&\pp_{\varepsilon-1}\cr0&0\cr\pp_{b+\varepsilon}&\pp_{b+\varepsilon-1}}$ and $W=YV$ or $YM_0V$.

If $\pp_0\ \pp_{a-b+1}\ge {\pp_a}^2$, the uniform convergence -- on ${{\cal A}'}^\NN$ -- of the sequence 
$\mathtt{N}(P_{\zeta_1}\dots P_{\zeta_k}YV)$ and $\mathtt{N}(P_{\zeta_1}\dots P_{\zeta_k}YM_0V)$ to the same
vector $V(\zeta)$ is insured by
Theorem \ref{convergence} and~Corollary~\ref{corol0}. When $n\to\infty$
the numerator in (\ref{over}) converges to
$V_1(\xi):=\pmatrix{1&1&0}M^*_{\xi_1}M_0^kQ_\varepsilon V(\zeta)$, and the denominator to $V_2(\xi):=\pmatrix{1&1&0}M_0^kQ_\varepsilon V(\zeta)$; this convergence is 
uniform on each cylinder $\[\varepsilon'0^k\varepsilon\]$. Since the first entrie in $Q_\varepsilon V(\zeta)$ is at least 
$\min\{\pp_\varepsilon,\pp_{\varepsilon-1}\}>0$, $V_1(\xi)$~and~$V_2(\xi)$ are positive and consequently $\phi_n(\xi)$ converges uniformly to 
$\log\displaystyle{V_1(\xi)\over V_2(\xi)}$.
This is also true on any finite reunion of such cylinders; let us denote by $X(k_0)$ the reunion of $\displaystyle\[\varepsilon'0^k\varepsilon\]$ for 
$k<k_0,\varepsilon\in{\cal A}\setminus\{0\}$ and
$\varepsilon'\in{\cal A}$; then
\begin{equation}\label{from1}
\forall\eta>0,\ \exists n_0\in\NN,\ 
n\ge n_0\hbox{ and }\xi\in X(k_0) \Rightarrow\left\vert\phi_n(\xi)-\log{V_1(\xi)\over V_2(\xi)}\right\vert\le\eta.
\end{equation}

We consider now a sequence $\xi\in{\cal A}^\NN\setminus X(k_0)$.
By using the left and right eigenvectors of~$M_0$ -- for the eigenvalue $\pp_0$ -- we obtain
$$
\lim_{k\to\infty}A_k=\lambda_0\pmatrix{\pp_0^2-\pp_a\pp_{b-1}&0&0\cr\pp_a\pp_b&0&0\cr\pp_0\pp_b&0&0}\quad\hbox{where }
\lambda_0>0,\ A_k:=\left\{\begin{array}{ll}\pp_0^{-k}{M_0}^k&\hbox{if }\pp_a\ \pp_{b-1}<{\pp_0}^2\\k^{-1}\pp_0^{-k}{M_0}^k&\hbox{if }\pp_a\ \pp_{b-1}={\pp_0}^2.
\end{array}\right.
$$
The entries $\pp_a\pp_b$ and $\pp_0\pp_b$ being positive, there exists $\lambda(\varepsilon')>0$ such that 
$$
\lim_{k\to\infty}\pmatrix{1&1&0}M^*_{\varepsilon'}A_kQ_\varepsilon
=\lambda(\varepsilon')\pmatrix{\pp_\varepsilon&\pp_{\varepsilon-1}}.
$$
Moreover the convergence of $\displaystyle\pmatrix{1&1&0}M^*_{\varepsilon'}A_kQ_\varepsilon\pmatrix{x\cr y}$ to 
$\lambda(\varepsilon')(\pp_\varepsilon x+\pp_{\varepsilon-1}y)$ is
uniform on the set of normalized nonnegative
column-vectors $\pmatrix{x\cr y}$. Similarily, there exists $\lambda_1>0$ such that
$\displaystyle\pmatrix{1&1&0}A_kQ_\varepsilon\pmatrix{x\cr y}$ converges uniformly to
$\lambda_1(\pp_\varepsilon x+\pp_{\varepsilon-1}y)$. Both limits are positive if $\varepsilon\ne0$, implying that the ratio converges uniformly to 
$\displaystyle{\lambda(\varepsilon')\over\lambda_1}$. Hence, using (\ref{over}) one can choose $k_0$ such that
-- if we assume $\xi\in{\cal A}^\NN\setminus X(k_0)$ and $\sigma\xi\ne\overline0$
\begin{equation}\label{from2}
n\ge k_0+4\Rightarrow\left\vert\phi_n(\xi)-\log{\lambda(\varepsilon')\over\lambda_1}\right\vert\le\eta.
\end{equation}
The uniform convergence of $\phi_n(\xi)$ on ${\cal A}^\NN$ follows from (\ref{from1}) and (\ref{from2}) since, in the remaining case $\sigma\xi=\overline0$ one has 
$\displaystyle\lim_{n\to\infty}\phi_n(\xi))=\log{\lambda(\xi_1)\over\lambda_1}$.


Conversely, suppose $\pp_a\ \pp_{b-1}>~{\pp_0}^2$. If $\mu_\pp^*$ is weak Gibbs
w.r.t. $\{\SS_\varepsilon\}_{\varepsilon=0}^{\ss-1}$ then, from~(\ref{A2})~and~(\ref{nstep}) one has $\phi_n(\xi)=o(n)$ for any $\xi\in{\cal A}^\NN$.
But this is not true: $\displaystyle\phi_{2n+1}(1\overline0)\sim n\log{\pp_a\pp_{b-1}\over{\pp_0}^2}$.

Suppose now $\pp_0\ \pp_{a-b+1}<{\pp_a}^2$. If $b=1$ we have $\pp_0<\pp_a$ hence $\pp_a\ \pp_{b-1}>~{\pp_0}^2$; that is, we are in the previous case. If $b\ne1$,
$\mu_\pp^*$ is no more weak Gibbs
w.r.t. $\{\SS_\varepsilon\}_{\varepsilon=0}^{\ss-1}$ because there exists a contradiction between the limit in (\ref{A2}) and the following:
$$
\lim_{n\to\infty}\left({\mu_\pp^*\[(0(a-b+1))^n1^n\]\over\mu_\pp^*\[(0(a-b+1))^n\]\cdot\mu_\pp^*\[1^n\]}\right)^{1/n}
={\pp_0\pp_{a-b+1}\over{\pp_a}^2}<1.
$$
~\qed

\end{section}

\eject
\baselineskip=14pt

\'Eric {\sc Olivier}\par
Centre de Ressources Informatiques\par
Universit\'e de Provence\par
3, place Victor Hugo\par
13331 MARSEILLE Cedex 3, France\par
{\sl E-mail : }Eric.Olivier@up.univ-mrs.fr

\medskip

Alain {\sc Thomas}\par
Centre de Math\'ematiques et d'Informatique\par
LATP \'Equipe de th\'eorie des nombres\par
39, rue F. Joliot-Curie\par
13453 Marseille Cedex 13, France\par
{\sl E-mail : }thomas@cmi.univ-mrs.fr

\end{document}